\documentclass[a4paper]{amsart} 
\usepackage{amssymb}
\usepackage{amscd}
\usepackage{times}

\numberwithin{equation}{subsection}

    \def\n{|\!|}    
\def\om{\omega}         \def\cb{\mathbb{C}}        
   \def\pb{\mathbb{P}}   \def\rb{\mathbb{R}}
\def\zbb{\mathbb{Z}}      \def\sc{{\mathcal{S}}}
\def\ot{\otimes}             
\def\zb{\overline{z}}     \def\phi{\varphi}     \def\spec{{Spec\    }}
\def\p1{{\pb^1_\zbb}}          \def\pbu{{\pb^1}}          \def\sn{S_n}
\def\scn{\mathcal{S}_n}                     \def\osn{\mathcal{O}_{S_n}}
\def\oscn{\mathcal{O}_{\mathcal{S}_n}}   \def\dega{{\widehat{deg}\  }}
\newtheorem{theo}{Theorem}               
           \newtheorem{lemma}{Lemma}

\begin{document}
\title{Analytic  torsion  of  Hirzebruch surfaces}  \author{Christophe
Mourougane} \address{Institut de Math{\'e}matiques de Jussieu / Plateau 7D
/ 175, rue du Chevaleret / 75013 Paris}
\email{christophe.mourougane@math.jussieu.fr} \maketitle

{\textbf{Abstract.}
\textit{ Using different forms of the arithmetic Riemann-Roch theorem
  and the computations of Bott-Chern secondary classes, we compute 
  the analytic torsion and the height of Hirzebruch surfaces.}
~

\section{Introduction}
In  ~\cite{RS},  Ray  and  Singer  defined the  analytic  torsion  for
 hermitian  complex manifolds,  which  can be  seen  as a  regularized
 determinant  of the  Laplacian for  the  $\overline\partial$ operator
 acting  on  the spaces  of  forms. They  computed  it  for curves  as
 functions  on   the  moduli  spaces  and  related   it  with  modular
 forms.  Yoshikawa generalized  this  relation for  Theta divisors  on
 Abelian varieties~\cite{yo}.
Another kind  of result is about Hermitian  symmetric manifolds, where
the spectrum of the Laplacian  (and hence the analytic torsion) can be
explicitly    derived    from   the    spectrum    of   the    Casimir
operator~\cite{ko}.

 We compute the analytic  torsion $\tau(S_n)=-\log T_0 (S_n,1)$ of (the
trivial flat  Hermitian line bundle on) the  Hirzebruch surfaces $S_n$
the only ruled surfaces over $\pbu$, endowed  with some canonically  
defined hermitian  metric.  We prove

\textbf{Main theorem}
$$
\tau(S_n)-\log Vol(S_n)=\frac{n\log(n+1)}{24}-\frac{n}{6}+2\tau(\pb^1).
$$

For  the analytic  structure  of  $S_n$ varies  with  $n$, results  on
adiabatic limits  (as in ~\cite{ber-bi}  for example) do not  apply at
once. There is nevertheless a striking resemblance with the asymptotic
of  the analytic  torsion associated  with the  positive  line bundles
$\mathcal{O}_\pbu (n)\to\pbu$ computed in~\cite{BV} Theorem 8.

Our  method  consists  in   applying different forms of the arithmetic
Riemann-Roch theorem, reducing the computation of the analytic torsion
to the computations of  distinguished Bott-Chern classes.  For we have
a description of the arithmetic Chow group of the natural $\zbb$-model
$\scn$  of the  Hirzebruch  surface  $S_n$, we  can  also compute  its
height.

\thanks{ I  would like to  thank Ken-Ichi Yoshikawa for  suggesting me
this   problem,  Daniel   Huybrechts, Xiaonan   Ma   and Vincent
Maillot    for   useful discussions.}


\section{Preliminaries}
\subsection{With the relative Euler sequence}
We  will recall  some  facts about  the  arithmetic Chow  ring of  the
projective line $\p1$ and of the Hirzebruch surfaces $\scn$.
We refer to the book~\cite{soule} for basics about the arithmetic
constructions. 

Let $\mathcal{V}$  be a  free $\zbb$-module of  rank $2$  and $\p1=\pb
(\mathcal{V})\to \spec\zbb$ the projective space of rank one quotients
of  $\mathcal{V}$.   The  choice  of  a Hermitian  scalar  product  on
$V:=\mathcal{V}\ot_\zbb \cb$ determines canonically a Hermitian metric
$h$ on  the tautological line bundle  $\mathcal{O}_{\pb^1}(1)$, as the
quotient metric  through $\pbu\times V\to\mathcal{O}_{\pb^1}(1)$.
Its curvature form  gives a K{\"a}hler form on $\pbu$  and will be denoted
by   $x=\omega_\pbu$.    The   arithmetic   first   Chern   class   of
$(\overline{\mathcal{O}_{\p1}(1)},h)$     will    be     denoted    by
$\widehat{x}$.        Considering       the       exact       sequence
$\overline{\mathcal{S}}$ with natural metrics
$$
0\to\overline{\mathcal{O}_{\p1}(-1)}\to \overline{\mathcal{V}^\star}
\to\overline{\mathcal{O}_{\p1}(1)}\to 0
$$   we  derive  the  relation
$(1-\widehat{x})(1+\widehat{x})-1=a(\widetilde{c}(\overline{\mathcal{S}}))
$    that    is,     using    ~\cite{gs-annals}    proposition~5.3,
$\widehat{x}^2=-a(\widetilde{c}_2(\overline{\mathcal{S}}))=a(x)$.

Consider      the     ample      vector      bundle     $\mathcal{E}_n
:=\mathcal{O}_{\p1}(1)\oplus \mathcal{O}_{\p1}(n+1)$  on $\p1$.  Then,
$$\mathcal{S}_n=\pb                    (\mathcal{E}_n)\stackrel{\pi}\to
\p1\stackrel{f}\to  \spec\zbb$$ is  a $\zbb$-model  of  the Hirzebruch
surface $S_n$.  The tautological  ample line bundle on $\mathcal{S}_n$
will be denoted  by $\mathcal{O}_{\mathcal{E}_n}(1)$.  From the metric
$h$,  we  construct  an  orthogonal  sum  Hermitian  metric  $h_n$  on
$E_n:=\mathcal{E}_n(\cb)$.    The    arithmetic    Chern   class    of
$\overline{\mathcal{E}_n}$                   is                  hence
$\widehat{c}(\overline{\mathcal{E}_n})=(1+\widehat{x})(1+(n+1)\widehat{x})$.
We also  get a Hermitian metric  (still denoted by  $h_n$) of positive
curvature        $\alpha_n:=\Theta(\mathcal{O}_{E_n}(1),h_n)$       on
$\mathcal{O}_{E_n}(1):=\mathcal{O}_{\mathcal{E}_n}(1)(\cb     )$    as
quotient  metric through $\pi^\star  E_n\to\mathcal{O}_{E_n}(1)$.  The
arithmetic           first           Chern          class           of
$\overline{\mathcal{O}_{\mathcal{E}_n}(1)}$   will   be   denoted   by
$\widehat{\alpha_n}$.   Considering   the   relative  Euler   sequence
$\overline{\Sigma(1)}$ metrized with  induced and quotient metric from
the metric $h_n$ on $E_n$
$$0\to      \overline{\mathcal{O}_{\pb      (\mathcal{E}_n)      }}\to
\pi^\star\overline{\mathcal{E}_n^\star}
\ot\overline{\mathcal{O}_{\mathcal{E}_n}(1)}\stackrel{q}\to
\overline{T_{\pb (\mathcal{E}_n)/\p1}}\to 0$$ we derive the relations
\begin{eqnarray}\label{c2trel}
\nonumber     0&=&\widehat{c_2}(\pi^\star\overline{\mathcal{E}_n^\star}
\ot\overline{\mathcal{O}_{\mathcal{E}_n}(1)})
+a(\widetilde{c_2}(\overline{\Sigma(1)}))\\
&=&\widehat{\alpha_n}^2-(n+2)\pi^\star
\widehat{x}\widehat{\alpha_n}+(n+1)\pi^\star\widehat{x}^2
+a(\widetilde{c_2}(\overline{\Sigma(1)})) .
\end{eqnarray}
\begin{eqnarray}
\label{c1trel}&&
\widehat{c_1}(\overline{T_{\scn/\p1}},\omega_q)
=\widehat{c_1}(\pi^\star\overline{\mathcal{E}_n^\star}
\ot\overline{\mathcal{O}_{\mathcal{E}_n}(1)})
+a(\widetilde{c_1}(\overline{\Sigma(1)}))
=2\widehat{\alpha_n}-(n+2)\pi^\star\widehat{x}.
\end{eqnarray}
Recall from our computations~\cite{ch}, that
$$\widetilde{c_2}(\overline{\Sigma(1)})=\widetilde{c_2}(\overline{\Sigma})=
-\Omega=           -\alpha_n           -\frac{\langle           \Theta
(E_n^\star,h)a^\star,a^\star\rangle_h}                         {\langle
a^\star,a^\star\rangle_h}$$   where  $\Omega$  denotes   the  relative
Fubini-Study form.  We have proved the following arithmetic analogs of
relations in Chow groups.
\begin{lemma}\label{lemme}
  \begin{list}{}{}
  \item   In   the   arithmetic   Chow   ring  of   $\p1$,   we   have
$$\widehat{x}^2=a(x).$$
\item In the arithmetic Chow ring of $\scn$, we have
$$\widehat{\alpha_n}^2-(n+2)\widehat{x}\widehat{\alpha_n}
=a\left(\alpha_n-(n+1)x              +\frac{\langle             \Theta
(E^\star,h)a^\star,a^\star\rangle_h}                           {\langle
a^\star,a^\star\rangle_h}\right).$$
  \end{list}
\end{lemma}

This exact sequence  also enables to compute the  height of Hirzebruch
surfaces       with       respect       to      the       polarization
$\overline{\mathcal{O}_{\mathcal{E}_n}(1)}$.  The arithmetic height of
an   arithmetic  variety  $\mathcal{X}\stackrel{p}\to   \spec\zbb$  of
relative    dimension   $n$   with    respect   to    a   polarization
$\overline{\mathcal{L}}\to\mathcal{X}$     is     defined    to     be
$\widehat{h}_{\overline{\mathcal{L}}}(\mathcal{X})
:=\dega\widehat{p}_\star\left(
\widehat{c_1}(\overline{\mathcal{L}})^{n+1}  \right)$  where $\dega  :
\widehat{CH}^1(\zbb)\to  \rb$  sends  $(\sum  n_\mathcal{P}\mathcal{P}
,\lambda  )$ to  $\sum n_\mathcal{P}\log\mathcal{P}  +\lambda/2$.  The
arithmetic height of $\scn $ is hence given by
$$\widehat{h}_{\overline{\mathcal{O}_{\mathcal{E}_n}(1)}}(\scn)
=\dega\widehat{f}_\star\widehat{\pi}_\star                       \left(
\widehat{\alpha_n}^3\right)                     =\dega\widehat{f}_\star
\widehat{s'_2}(\overline{\mathcal{E}_n})$$                        where
$\widehat{s'_m}(\overline{\mathcal{E}_n}):=\widehat{\pi}_\star   \left(
\widehat{\alpha_n} ^{1+m} \right)$ is the $m$-th geometric Segre class
of   $\overline{\mathcal{E}_n}$   in   the   arithmetic   Chow   group
$\widehat{CH}^m(\p1)$ of $\p1$.  Define the secondary form
$$S_{m+1}(E_n,h_n)
=\pi_\star(\widetilde{c}_2(\overline{\Sigma}(1))\alpha_n^m)
\in\widetilde{A^{m,m}(\pb^1)}.
$$ From  the degree two  relation ~(\ref{c2trel}) on  arithmetic Chern
classes,  taking  product with  $\widehat{\alpha_n}^m$  ($m=0$ ;  $1$)
pushing-forward through $\widehat{\pi}_\star$ and applying usual rules
of  calculus on  arithmetic Chow  groups we  are led  to  relations in
$\widehat{CH}(\p1 )$
 \begin{eqnarray*}
\widehat{s'_1}(\overline{\mathcal{E}_n})-(n+2)\widehat{x}+a(S_{1}(E_n,h_n))
&=&0\\    \widehat{s'_2}(\overline{\mathcal{E}_n})   -(n+2)\widehat{x}
\widehat{s'_1}(\overline{\mathcal{E}_n})                              +
(n+1)\widehat{x}^2+a(S_{2}(E_n,h_n))&=&0
 \end{eqnarray*}
and then to
\begin{eqnarray*}\label{eq}
\widehat{s'_2}(\overline{\mathcal{E}_n})       =(n^2+3n+3)\widehat{x}^2
-(n+2)a(xS_1(E_n,h_n))-a(S_{2}(E_n,h_n)).
\end{eqnarray*}
Using,       $\widetilde{c}_2(\overline{\Sigma}(1))=-\Omega$       and
$\displaystyle\int_\cb\frac{idz\wedge d\zb}{2\pi(1+|z|^2)^3}=\frac{1}{2}$ 
we infer,
  \begin{eqnarray*}
S_1(E_n,h_n)&=&\pi_\star(\widetilde{c}_2(\overline{\Sigma}(1)))
=-\pi_\star(\Omega)=                                               -1\\
S_{2}(E_n,h_n)&=&-\pi_\star(\Omega\alpha_n)=\pi_\star(\Omega
\frac{\langle  \Theta  (E_n^\star,h)a^\star,a^\star\rangle_h} {\langle
a^\star,a^\star\rangle_h})=\frac{1}{2}c_1(E^\star_n,h_n).
 \end{eqnarray*}
Hence,
$$\widehat{s'_2}(\overline{\mathcal{E}_d})=(n^2+3n+3)\widehat{x}^2
+\frac{3(n+2)}{2}a(x)=\frac{2n^2+9n+12}{2}a(x).$$  Finally,  computing
the arithmetic degree we get
\begin{theo} The arithmetic height of the Hirzebruch surface $\scn$ is
\begin{eqnarray*}
\widehat{h}_{\overline{\mathcal{O}_{\mathcal{E}_n}(1)}}(\scn)
=\frac{2n^2+9n+12}{4}.
 \end{eqnarray*}
\end{theo}

\subsection{With the sequence associated to the fibration}
To  complete the  picture of  the arithmetic  datas of  the Hirzebruch
surfaces, we need  to compute the Bott-Chern secondary  classes of the
short exact sequence
$$0  \to  (T_{S_n/\pbu},\alpha_n) \stackrel{\iota}\to  (TS_n,\alpha_n)
\stackrel{d\pi}\to  (\pi^\star T\pbu,\alpha_n)  \to  0$$ of  hermitian
vector bundles on $S_n$.

We       need       to       make      explicit       the       metric
$\alpha_n=\Theta(\mathcal{O}_{E_n}(1),h_n)$.   We   choose   a   local
holomorphic  frame  $e$ for  $\mathcal{O}_{E_n}(1)$.  This provides  a
frame  $(e^\star,{e^\star}^{n+1})$  for   $E^\star_n$.  We  denote  by
$(a_1,a_2)$ the  corresponding coordinates  on $E^\star_n$ and  by $z$
the holomorphic coordinate $z:=a_1/a_2$  on an appropriate open set of
$\pb   (E_n)$.  Then,   $\alpha_n$  is   computed  by   $dd^c\log  (\n
e^\star+z{e^\star}^{n+1}\n^2_{h_n})$.  We find
\begin{eqnarray}
\label{alpha} \nonumber
\alpha_n&=&\frac{1+(n+1)|z|^2\n       e^\star\n^{2n}}       {1+|z|^2\n
e^\star\n^{2n}}dd^c\log\n     e^\star      \n^2\\     &&     +\frac{\n
e^\star\n^{2n}}{(1+|z|^2\n                           e^\star\n^{2n})^2}
\frac{i}{2\pi}\Big|dz+nzd'\log\n e^\star \n^2 \Big|^2.
\end{eqnarray}
The choice  of the frame  $e$ naturally induces  a choice for  a local
holomorphic splitting  $TS_n\simeq \pi^\star T\pbu\oplus T_{S_n/\pbu}$
so   that   for  a   vector   $X$   in   $T\pbu$,  $d\pi(X,0)=X$   and
$i\left(\frac{\partial}{\partial z}\right)=(0,\frac{\partial}
{\partial z})$. We find
\begin{eqnarray*}
\n     (X,0)\n^2_{\alpha_n}&=&\frac{1+(n+1)|z|^2\n     e^\star\n^{2n}}
{1+|z|^2\n   e^\star\n^{2n}}\n   X\n^2_{\pbu}  +\frac{1}{2\pi}\frac{\n
e^\star\n^{2n}}{(1+|z|^2\n       e^\star\n^{2n})^2}      |nzd'_X\log\n
e^\star\n^{2}|^2\\      \langle      (X,0),(0,\frac{\partial}{\partial
z})\rangle  _{\alpha_n}   &=&\frac{1}{2\pi}\frac{  \n  e^\star\n^{2n}}
{(1+|z|^2\n    e^\star\n^{2n})^2}nzd'_X   \log\n    e^\star\n^2\\   \n
(0,\frac{\partial}{\partial             z})\n^2            _{\alpha_n}
&=&\frac{1}{2\pi}\frac{\n                    e^\star\n^{2n}}{(1+|z|^2\n
e^\star\n^{2n})^2}   =\frac{1}{2\pi}\frac{\n  e^\star\n^{2(n+2)}}  {\n
e^\star + z {e^\star}^{n+1}\n ^4}.
\end{eqnarray*}

Let $X$  be a vector in $T\pbu$.  Its image $(d\pi)^\star\pi^\star(X)$
can      be       written      as      $(d\pi)^\star\pi^\star      (X)
=(X,0)+V_X(0,\frac{\partial}{\partial  z})$ subject  to  the condition
$\langle   X+V_X\frac{\partial}{\partial   z},\frac{\partial}{\partial
z}\rangle   _{\alpha_n}=0$.   This   leads   to  $V_X=-nz   d'_X\log\n
e^\star\n^2$.   We can  then derive  the expression  for  the quotient
metric on $\pi^\star T\pbu$
\begin{eqnarray*}
\frac{{\alpha_n}_{|\pi^\star T\pbu}}{\pi^\star\omega_\pbu} &=&\frac{\n
(d\pi)^\star  X \n^2_{\alpha_n}}{\n  X\n^2_\pbu} =\frac{1+(n+1)|z|^2\n
e^\star\n^{2n}}{1+|z|^2\n        e^\star\n^{2n}}        =\frac{\langle
\pi^\star\Theta(E^\star_n,h_n)a^\star,a^\star\rangle}        {\pi^\star
\omega_\pbu \langle a^\star,a^\star\rangle}.
\end{eqnarray*}

We can now  compute the Bott-Chern class for  the short exact sequence
constructed from~$d\pi$.
\begin{lemma}\label{lem:BC}
The Bott-Chern class for Chern classes in the short exact sequence
$$0  \to  (T_{S_n/\pbu},\alpha_n) \stackrel{\iota}\to  (TS_n,\alpha_n)
\stackrel{d\pi}\to    (\pi^\star    T\pbu,\alpha_n)    \to   0$$    is
\begin{eqnarray*}           \widetilde{c}(TS_n,T\pbu,\alpha_n,\alpha_n)
&=&\left(  \frac{n}{1+(n+1)|z|^2\n e^\star\n^{2n}} -\frac{n}{1+|z|^2\n
e^\star\n^{2n}}\right)\pi^\star x.
\end{eqnarray*}
\end{lemma}
\begin{proof}
As   a   consequence  of~\cite{gs-annals}   prop   1.2.5,  we   derive
$\widetilde{c}_1(TS_n,T\pbu,\alpha_n,\alpha_n)=0$.  Using the original
computations of Bott and Chern,
\begin{eqnarray*}
\widetilde{c}_2(TS_n,T\pbu,\alpha_n,\alpha_n)                 =\int_0^1
\frac{\Phi(u)-\Phi(0)}{u}du
\end{eqnarray*}
with
\begin{eqnarray*}
\Phi(u)&=&Trace\left(\left(1-u\right)\Theta (\pi^\star T\pbu,\alpha_n)
+ u(d\pi)\Theta (TS_n,\alpha_n) (d\pi)^\star \right)\\ &=&(d\pi)\Theta
(TS_n,\alpha_n)    (d\pi)^\star\\&&    -\left(1-u\right)\frac{i}{2\pi}
(d\pi)(\nabla'_{Hom(T_{S_n/\pbu},TS_n)}\iota)               \wedge\iota
^\star\nabla''_{Hom(\pi^\star T\pbu,TS_n)} (d\pi)^\star,
\end{eqnarray*}
follows
\begin{eqnarray*}
\widetilde{c}_2(TS_n,T\pbu,\alpha_n,\alpha_n)         &=&\frac{i}{2\pi}
d\pi(\nabla'_{Hom(T_{S_n/\pbu},TS_n)}\iota)                 \wedge\iota
^\star\nabla''_{Hom(\pi^\star T\pbu,TS_n)} (d\pi)^\star.
\end{eqnarray*}
 Now, remark that
\begin{eqnarray*}
\iota  ^\star  (\nabla''   (d\pi)^\star)X  &=&\iota  ^\star  \nabla''(
(d\pi)^\star  X)  =d''V_X\iota  ^\star(0,\frac{\partial}{\partial  z})
=d''V_X \frac{\partial}{\partial z}.
\end{eqnarray*}
Also,
\begin{eqnarray*}
 d\pi(\nabla'\iota)            (\frac{\partial}{\partial           z})
&=&d\pi\nabla'\frac{\partial}{\partial               z}              =
\frac{\{\nabla'\frac{\partial}{\partial
z},X+V_X\frac{\partial}{\partial                    z}\}_{\alpha_n}}{\n
X+V_X\frac{\partial}{\partial                        z}\n^2_{\alpha_n}}
d\pi\left((X,0)+V_X(0,\frac{\partial}{\partial             z})\right)\\
&=&\frac{d'\langle\frac{\partial}{\partial                          z},
X+V_X\frac{\partial}{\partial    z}\rangle_{\alpha_n}-\{\frac{\partial}
{\partial   z},d''V_X\frac{\partial}{\partial   z}\}_{\alpha_n}}   {\n
X+V_X\frac{\partial}{\partial                       z}\n^2_{\alpha_n}}X
=-d'\overline{V_X}\frac{\n       \frac{\partial}{\partial       z}\n^2
_{\alpha_n}}{\n   X+V_X\frac{\partial}{\partial  z}\n^2_{\alpha_n}}X\\
&=&-\frac{1}{2\pi}     \frac{\n     e^\star\n^{2(n+2)}}    {(1+|z|^2\n
e^\star\n^{2n})(1+(n+1)|z|^2\n e^\star\n^{2n})} d'\overline{V_X}.
\end{eqnarray*}
Summing up
\begin{eqnarray*}
\lefteqn{\widetilde{c}_2(TS_n,T\pbu,\alpha_n,\alpha_n)
=-\frac{i}{4\pi^2}\frac{\n  e^\star\n^{2n}}{(1+|z|^2\n e^\star\n^{2n})
(1+(n+1)|z|^2\n  e^\star\n^{2n})}\frac{d'\overline{V_X}\wedge  d''V_X}
{\n     X\n^2_\pbu     }}&&\\     &=&-\frac{i}{4\pi^2}\frac{n^2|z|^2\n
e^\star\n^{2n}}{(1+|z|^2\n       e^\star\n^{2n})       (1+(n+1)|z|^2\n
e^\star\n^{2n})} \frac{d'd''_X\log\n e^\star\n^{2}\wedge d''d'_X\log\n
e^\star\n^2}        {\n       X\n^2_\pbu}\\       &=&-\frac{n^2|z|^2\n
e^\star\n^{2n}}{(1+|z|^2\n       e^\star\n^{2n})       (1+(n+1)|z|^2\n
e^\star\n^{2n})}\pi^\star\omega_\pbu.
\end{eqnarray*}
The last  equality is derived  from the fact that  those $(1,1)$-forms
take  the same  values  on $\pi^\star  X$  which generates  $\pi^\star
T\pbu$.
\end{proof}

\subsection{Some characteristic classes}
In the previous two exact sequences the bundle $T_{S_n/\pbu}$ appeared
equipped  with two  different metrics,  the one  $\omega_q$  gotten by
quotient of that on $ \pi^\star E^\star_n\ot \mathcal{O}_{E_n}(1)$ and
the  one induced  by $\alpha_n$.   We will  compare those  metrics and
derive expressions for some characteristic classes.

Choose  a  point   $x_0$  in  $\pbu$  and  a   normal  frame  $e$  for
$\mathcal{O}_{\pbu}(1)$ at $x_0$. Then, on the fiber of~$x_0$
 $${\alpha_n}_{|T_{S_n/\pbu}}         =        d_zd_z^c\log        (\n
e^\star+z{e^\star}^{n+1}\n^2)=         d_zd_z^c\log         (1+|z|^2)=
\frac{i}{2\pi}\frac{dz\wedge  d\zb}{(1+|z|^2)^2}.$$ On the  other end,
the  map  $q$ is  built  from the  differential  of  the quotient  map
$E_n^\star-\pbu\times\{0\}           \to           \pb          (E_n),
(x,a_1e^\star+a_2{e^\star}^{n+1})\mapsto      (x,[a_1:a_2])     \simeq
(x,z:=\frac{a_2}{a_1})$
$$\begin{array}{cccc}    q:&   E_n^\star\ot\mathcal{O}_{E_n}(1)   &\to
&T_{\pb                         (E_n)/\pbu}\\                        &
(b_1e^\star+b_2{e^\star}^{n+1})\ot(a_1e^\star+a_2{e^\star}^{n+1})^\star
&\mapsto&\frac{a_1b_2-a_2b_1}{a_1^2}\frac{\partial}{\partial z}.
\end{array}$$
Follows the expression for its adjoint map
\begin{eqnarray*}
q^\star                     \left(\frac{1}{a_1}\frac{\partial}{\partial
    z}\right)&=&\left({e^\star}^{n+1}                    -\frac{\langle
    {e^\star}^{n+1},a_1e^\star+a_2{e^\star}^{n+1}\rangle     }     {\n
    a_1e^\star+a_2{e^\star}^{n+1}\n^2}(a_1e^\star+a_2{e^\star}^{n+1})\right)
    \otimes     (a_1e^\star+a_2{e^\star}^{n+1})^\star    \\    q^\star
    \left(\frac{\partial}{\partial  z}\right)  &=&\frac{-\zb  e^\star+
    {e^\star}^{n+1}}{1+|z|^2} (e^\star+ z{e^\star}^{n+1})^\star.
\end{eqnarray*}
We    derive    that    the    quotient    metric    is    given    by
  $\displaystyle\omega_q=i\frac{dz\wedge    d\zb}{(1+|z|^2)^2}   =2\pi
  {\alpha_n}_{|T_{S_n/\pbu}}$ on the fiber  of $x_0$.  This enables to
  express     the     arithmetic      first     Chern     class     of
  $\overline{T_{\sc_n/\p1}}$          with          the         metric
  ${\alpha_n}_{|T_{S_n/\pbu}}$
\begin{eqnarray*}
  \widehat{c_1}(\overline{T_{\scn/\p1}},\alpha_n)
&=&\widehat{c_1}(\overline{T_{\scn/\p1}},\omega_q)+a(\log       2\pi)\\
&=&2\widehat{\alpha_n}-(n+2)\widehat{x} +a(\log 2\pi)
\end{eqnarray*}
where we  have used  the relation~(\ref{c1trel}).  Similarly,  the map
$(T\pbu,2\pi\omega_\pbu)\to     (\mathcal{O}_{\pbu}(2),h)$    is    an
isometry. Just recall
$$                               \widehat{c}(\overline{T\scn},\alpha_n)
=\widehat{c}(\overline{T_{\scn/\p1}},\alpha_n)
\widehat{c}(\overline{T_\p1                                 },\alpha_n)
-a(\widetilde{c}(TS_n,T\pbu,\alpha_n,\alpha_n))$$ to end the proof of
\begin{lemma}\label{lem:classes}
\begin{eqnarray*}
  \widehat{c_1}(\overline{T_{\scn/\p1}},\alpha_n)
&=&2\widehat{\alpha_n}-(n+2)\widehat{x}         +a(\log        2\pi)\\
\widehat{c_1}(\overline{T\p1         },\alpha_n)&=&2\widehat{x}-a\left(
\log\frac{{\alpha_n}_{|\pi^\star                  T\pbu}}{2\pi\pi^\star
\omega_\pbu}\right)\\          \widehat{c_1}(\overline{T\scn},\alpha_n)
&=&2\widehat{\alpha_n}-n\widehat{x}                            -a\left(
\log\frac{{\alpha_n}_{|\pi^\star                T\pbu}}{4\pi^2\pi^\star
\omega_\pbu}\right)\\          \widehat{c_2}(\overline{T\scn},\alpha_n)
&=&4\widehat{\alpha_n}\widehat{x}-2(n+2)\widehat{x}^2+a\Big(2\log 2\pi
x\Big.\\           &&\Big.           -\log\frac{{\alpha_n}_{|\pi^\star
T\pbu}}{2\pi\pi^\star                                      \omega_\pbu}
c_1(T_{S_n/\pbu})-\widetilde{c}_2(TS_n,T\pbu,\alpha_n,\alpha_n)\Big).
\end{eqnarray*}
\end{lemma}

\section{Using the arithmetic Riemann-Roch theorem}

We    will   apply    the   degree    one    arithmetic   Riemann-Roch
 theorem~\cite{gs-invent} for the map $F : \mathcal{S}_n\to \spec\zbb$
 to  compute  the  analytic  torsion  of the  bundles  of  holomorphic
 differential forms $\Omega^p_{S_n}$ on $S_n$.

\subsection{On the analytic torsion}
We recall the definition of the analytic torsion of a Hermitian
vector bundle $(E,h)$ on a compact K{\"a}hler manifold $(X,\omega)$.
We endow the space of differential forms with values in $E$ with its
$L^2$ metric (constructed with $\omega$ and $h$) in order to construct
the adjoint $\overline{\partial_q}^\star$ of the Dolbeault operator
$\overline{\partial_q}~:~A^{0,q} (X,E)\to A^{0,q+1} (X,E)$. The
associated Laplace operator $\Delta_q''$ has a discrete spectrum 
$$0=0=\cdots=0<\lambda_1\leq\lambda_2\leq\cdots\leq\lambda_N\leq\cdots.$$  
Its spectral function $\displaystyle
\zeta_q(s):=\sum_1^{+\infty}\lambda_N^{-s}$  extends to a meromorphic
function on $\cb$, holomorphic at $0$. We define the regularized
determinant of the restriction of the Laplace operator to the
orthogonal complement of its kernel to be
$\displaystyle
det' \Delta_q'':=\exp\left(-\frac{d\ }{ds}\zeta_q(0)\right)$.
The analytic torsion of $(X,\omega ,E,h)$ is then given by
\begin{eqnarray*}
T_0(X,\omega ,E,h):=\exp(-\tau((X,\omega ,E,h))
:=\prod_{q\geq 0}(det' \Delta_q'')^{(-1)^q q}.
\end{eqnarray*}

\subsection{On arithmetic first Chern class of determinant bundles}
On one  hand, from the very  definition of the  arithmetic first Chern
class  of a  Hermitian  line bundle  and  from the  definition of  the
Quillen metric on  the determinant of the direct  image of a Hermitian
vector bundle,  identifying $CH^1(\spec\zbb )$ with  $\rb$ through the
degree map $\dega$, we get~(see~\cite{gs-invent} 4.1.5)
\begin{eqnarray*}
\lefteqn{\dega                 \widehat{c_1}(\lambda_F(\overline{\Omega
_{\scn}^p}),Quillen    )}\\   &=&\sum_{q=0}^2   (-1)^q\left(\log\sharp
H^q(S_n,\Omega_{\sn}                ^p)_{torsion}                -\log
Vol_{L^2}\frac{H^q(S_n,\Omega_{\sn}     ^p)}     {H^q(S_n,\Omega_{\sn}
^p)_\zbb} \right) +\frac{1}{2}\tau (S_n,\Omega_{\sn} ^p).
\end{eqnarray*}

From Leray-Hirsch theorem, the cohomology ring $H^{\bullet} (S_n ,\zbb
)$ of $S_n=\pb (E_n)$ is
$$\frac{H^{\bullet}(\pbu             ,\zbb            )[\{\alpha_n\}]}
{\{\alpha_n\}^2-(n+2)\{\pi^\star\omega_\pbu\}\cup\{\alpha_n\}}
=\frac{\zbb[\{\pi^\star\omega_\pbu\},\{\alpha_n\}]}
{\{\pi^\star\omega_\pbu\}^2,
\{\alpha_n\}^2-(n+2)\{\pi^\star\omega_\pbu\}\cup\{\alpha_n\}}.
$$  We  hence  get,  applying  Hodge  decomposition  to  the  De  Rham
cohomology,
\begin{list}{}{}
\item  for  $p=0$, $H^{\bullet}  (S_n  ,\osn  )_\zbb=H^{0} (S_n  ,\osn
)_\zbb=\zbb 1$
\item       for      $p=1$,      $H^{\bullet}(S_n,\Omega_{\sn}^1)_\zbb
=H^1(S_n,\Omega_{\sn}
^1)_\zbb=\zbb\{\pi^\star\omega_\pbu\}+\zbb\{\alpha_n\}$
\item     for     $p=2$,    $H^{\bullet}(S_n,\Omega_{\sn}     ^2)_\zbb
=H^2(S_n,\Omega_{\sn}
^2)_\zbb=\zbb\{\pi^\star\omega_\pbu\}\cup\{\alpha_n\}
=\zbb\frac{\{\alpha_n\}^2}{n+2} $
\end{list}
We  now   intend  to  find  the  harmonic   representatives  of  those
generators,  with respect  to the  metric $\alpha_n$.  The  forms $1$,
$\alpha_n$ and $\frac{\alpha_n^2}{n+2}$ are easily seen to be harmonic
of       norm      $\sqrt{\frac{n+2}{2}}$,       $\sqrt{n+2}$      and
$\sqrt{\frac{2}{n+2}}$  respectively.   Computation  of  the  harmonic
representative  of $\{\pi^\star\omega_\pbu\}$  is  quite involved  and
will not  actually be  used.  Details will  be given in  the appendix.
The                form                $\omega_H:=\pi^\star\omega_\pbu
-\frac{1}{n+2}dd^c\log\frac{\langle                              \Theta
(E^\star,h)a^\star,a^\star\rangle_h}       {\pi^\star\omega_\pbu\langle
a^\star,a^\star\rangle_h}$ is harmonic.  To compute its $L^2$-norm, we
need to introduce  the Hodge $\star$-operator. It is  defined in order
to    fulfill    the    relation    $u\wedge    \star    v=    \langle
u,v\rangle_{\alpha_n}dV_{\alpha_n}$.             Recall           that
$dV_{\alpha_n}=\alpha_n^2/2$         and        that        $|\alpha_n
|^2_{\alpha_n}=2$.  Local computation  of  the Hodge  $\star$-operator
leads  to  $\star  \alpha_n=\alpha_n$.   For $\star\omega_H$  is  also
harmonic,  it can  be  written as  $a\alpha_n+b\omega_H$. Writing  the
relation $\star\star\omega_H=\omega_H$  leads to either  $(b=1 ; a=0)$
or     $b=-1$.      Noting     that     $b=1$    would     lead     to
$\n\omega_H\n^2=\int\omega_H           \wedge\star\omega_H=\int\omega_H
\wedge\omega_H   =\int  \pi^\star\omega_\pbu\wedge\pi^\star\omega_\pbu
=0$,  we  derive  $b=-1$. Equating  $\omega_H\wedge\star\alpha_n$  and
$\star\omega_H\wedge\alpha_n$                  we                  get
$a\alpha_n^2=2\alpha_n\wedge\omega_H$.     Integrating,     we    find
$a=\frac{2}{n+2}$. We can now conclude
$$\n\omega_H\n^2=\int_{S_n}\omega_H\wedge(\frac{2}{n+2}\alpha_n-\omega_H)
=\frac{2}{n+2}.$$  We hence  find an  orthonormal basis  of $H^\bullet
(S_n,\Omega_{\sn}^1)      =\cb\{\alpha_n\}\oplus      H_{prim}^\bullet
(S_n,\Omega_{\sn}^1)$
$$\frac{\{\alpha_n\}}{\sqrt{n+2}}            \           ;           \
\sqrt{n+2}\left(\{\pi^\star\omega_\pbu\}-\frac{\{\alpha_n\}}{{n+2}}\right)$$
The      $L^2$-volume     of      $\frac{H^1(S_n,\Omega_{\sn}     ^1)}
{H^1(S_n,\Omega_{\sn} ^1)_\zbb}$,  which is the norm  of the generator
$\{\pi^\star\omega_\pbu\}\wedge\{\alpha_n\}$        of        $\Lambda
^2H^1(S_n,\Omega_{\sn} ^1)_\zbb$  is $1$. Concluding  this first step,
we get
\begin{eqnarray*}
\widehat{c_1}(\lambda_F(\overline{\Omega    _{\scn}^0}),   Quillen   )
&=&-\log\sqrt{\frac{n+2}{2}}+\frac{1}{2}\tau(\sn          )         \\
\widehat{c_1}(\lambda_F(\overline{\Omega    _{\scn}^1}),   Quillen   )
&=&-\frac{1}{2}\tau(\sn,\Omega          _{\sn}^1          )         \\
\widehat{c_1}(\lambda_F(\overline{\Omega    _{\scn}^2}),   Quillen   )
&=&\log\sqrt{\frac{n+2}{2}}+\frac{1}{2}\tau(\sn,\Omega _{\sn}^2 )
\end{eqnarray*}

\subsection{The arithmetic Riemann-Roch theorem}\label{RR}
On the other hand, the arithmetic Riemann-Roch theorem reads
\begin{eqnarray*}
\widehat{c_1}(\lambda_F(\overline{\Omega     _{\scn}^p},Quillen     ))
&=&\widehat{F}_\star\left(    \widehat{Td}^R(\overline{T\scn},\alpha_n)
\widehat{ch}(\overline{\Omega      _{\scn}^p},\alpha_n)      \right)\\
&=&\widehat{F}_\star\left(      \widehat{Td}(\overline{T\scn},\alpha_n)
\widehat{ch}(\overline{\Omega                      _{\scn}^p},\alpha_n)
\right)-a(F_\star\left( Td({TS_n})R({TS_n}) ch(\Omega _{S_n}^p)\right)
\end{eqnarray*}

Recall
\begin{eqnarray*}
\widehat{Td}(\overline{T\scn})=1+\frac{\widehat{c_1}(\overline{T\scn})}{2}
+\frac{\widehat{c_1}(\overline{T\scn})^2+\widehat{c_2}(\overline{T\scn})}{12}
+\frac{\widehat{c_1}(\overline{T\scn})\widehat{c_2}(\overline{T\scn})}{24}
\end{eqnarray*}
and
\begin{eqnarray*}
\widehat{ch}(\overline{\Omega                      _{\scn}^1},\alpha_n)
&=&2-\widehat{c_1}(\overline{T\scn})
+\frac{\widehat{c_1}(\overline{T\scn})^2-2\widehat{c_2}(\overline{T\scn})}{2}
-\frac{\widehat{c_1}(\overline{T\scn})^3
-3\widehat{c_1}(\overline{T\scn})\widehat{c_2}(\overline{T\scn})}{6}\\
\widehat{ch}(\overline{\Omega                      _{\scn}^2},\alpha_n)
&=&1-\widehat{c_1}(\overline{T\scn})
+\frac{\widehat{c_1}(\overline{T\scn})^2}{2}
-\frac{\widehat{c_1}(\overline{T\scn})^3}{6}
\end{eqnarray*}
so that
\begin{eqnarray*}
\left[\widehat{Td}(\overline{T\scn})
\widehat{ch}(\overline{\mathcal{O}_\p1},1)                    \right]_3
=-\left[\widehat{Td}(\overline{T\scn})\widehat{ch}(\overline{\Omega
_{\scn}^2},\alpha_n)                                          \right]_3
&=&\frac{\widehat{c_1}(\overline{T\scn})\widehat{c_2}(\overline{T\scn})}{24}\\
\left[\widehat{Td}(\overline{T\scn})\widehat{ch}(\overline{\Omega
_{\scn}^1},\alpha_n) \right]_3&=&0.
\end{eqnarray*}
From    formulas    for    the    arithmetic    Chern    classes    of
$(\overline{T\scn},\alpha_n)$ in lemma~\ref{lem:classes} we get
\begin{eqnarray*}
 {\widehat{c_1}(\overline{T\scn},\alpha_n)
\widehat{c_2}(\overline{T\scn},\alpha_n)}
&=&8\widehat{\alpha_n}^2\widehat{x}-8(n+1)\widehat{\alpha_n}\widehat{x}^2
+a\left( -  4\log\frac{{\alpha_n}_{|\pi^\star T\pbu}} {4\pi^2\pi^\star
\omega_\pbu}\alpha_n  x  +2\log  2\pi  xc_1(TS_n)  \right.\\  &&\left.
-\log\frac{{\alpha_n}_{|\pi^\star  T\pbu}}{2\pi\pi^\star  \omega_\pbu}
c_1(TS_n,\alpha_n)c_1(T_{S_n/\pbu})
-c_1(TS_n,\alpha_n)\widetilde{c}_2(TS_n,T\pbu,\alpha_n,\alpha_n\right).
\end{eqnarray*}

To compute the integral, recall from lemma~\ref{lemme} that
\begin{eqnarray*}
\widehat{F}_\star(\widehat{\alpha_n}\widehat{x}^2)&=&
a(F_\star(\alpha_n             x))=a(1)\\            \widehat{F}_\star
(\widehat{\alpha_n}^2\widehat{x})             &=&(n+2)\widehat{F}_\star
(\widehat{\alpha_n}\widehat{x}^2) +a(F_\star (\alpha_n x))=a(n+3)
\end{eqnarray*}
Applying  the  morphism $\omega$  in  lemma~\ref{lem:classes}, we  can
compute the curvature form of $(TS_n,\alpha_n)$.
\begin{eqnarray*}
\nonumber              c_1(TS_n,\alpha_n)              &=&2\alpha_n-nx
 -dd^c\log\frac{{\alpha_n}_{|\pi^\star    T\pbu}}{\pi^\star\omega_\pbu}
 \\   \nonumber   &=&\left(  n+2-\frac{3n}{1+|z|^2\n   e^\star\n^{2n}}
 +\frac{n}{1+(n+1)|z|^2\n      e^\star\n^{2n}}\right)\pi^\star     x\\
 &&+\left(\frac{3\n    e^\star\n^{2n}}{(1+|z|^2\n   e^\star\n^{2n})^2}
 -\frac{(n+1)\n                         e^\star\n^{2n}}{(1+(n+1)|z|^2\n
 e^\star\n^{2n})^2}\right)   \frac{i}{2\pi}\Big|dz+nzd'\log\n  e^\star
 \n^2 \Big|^2
\end{eqnarray*}
where  we used  the  expression~(\ref{alpha}) for  $\alpha_n$ and  the
expression  in the appendix  for $dd^c\log\frac{{\alpha_n}_{|\pi^\star
T\pbu}}{\pi^\star\omega_\pbu}$.       Then     note      that     from
lemma~\ref{lem:classes} we find
\begin{eqnarray*}
 c_1(T_{S_n/\pbu},\alpha_n)       &=&2\alpha_n-(n+2)\pi^\star      x\\
&=&(n-\frac{2n}{1+|z|^2\n   e^\star\n^{2n}})\pi^\star   x   +\frac{2\n
e^\star\n^{2n}}{(1+|z|^2\n    e^\star\n^{2n})^2}    \Big|dz+nzd'\log\n
e^\star \n^2 \Big|^2.
\end{eqnarray*}
so that
\begin{eqnarray*}
\frac{c_1(TS_n,\alpha_n)   c_1(T_{S_n/\pbu},\alpha_n)}   {\pi^\star  x
\wedge\Big|dz+nzd'\log\n   e^\star  \n^2   \Big|^2}  &=&\frac{(5n+2)\n
e^\star\n^{2n}}{(1+|z|^2\n    e^\star\n^{2n})^2}   +\frac{(n+1)(n+2)\n
e^\star\n^{2n}}{(1+(n+1)|z|^2\n     e^\star\n^{2n})^2}    -\frac{12n\n
e^\star\n^{2n}}{(1+|z|^2\n e^\star\n^{2n})^3}.
\end{eqnarray*}
Integration leads to
\begin{eqnarray*}\label{c1c1}
\lefteqn{F_\star  \left( c_1(TS_n,\alpha_n) c_1(T_{S_n/\pbu},\alpha_n)
\log\frac{{\alpha_n}_{|\pi^\star          T\pbu}}{\pi^\star\omega_\pbu}
\right)}&&\nonumber\\                                         \nonumber
&=&(5n+2)\left[\frac{n+1}{n}\log(n+1)-1 \right] \nonumber +(n+2)\left[
1-\frac{1}{n}\log(n+1)\right]\\                               \nonumber
&&-6(n+1)\left[\frac{n+1}{n}\log(n+1)-1           \right]          +3n
=5n+6-\left[n+6+\frac{6}{n} \right]\log(n+1)
\end{eqnarray*}
Thanks to lemma~\ref{lem:BC} we get
\begin{eqnarray*}\label{c1c2tilde}
\lefteqn{F_\star\left(               c_1(TS_n               ,\alpha_n)
\widetilde{c_2}(TS_n,T\pbu,\alpha_n,\alpha_n)\right)}     \nonumber&&\\
\nonumber                                                   &=&\int_\cb
\left(\frac{3}{(1+|z|^2)^2}-\frac{n+1}{(1+(n+1)|z|^2)^2}\right) \left(
\frac{n}{1+(n+1)|z|^2}-     \frac{n}{1+|z|^2}\right)dz\wedge    d\zb\\
\nonumber                                                   &=&\int_\cb
-\frac{3n}{(1+|z|^2)^3}-n\frac{n+1}{(1+(n+1)|z|^2)^3}
+(n+1)\frac{n+1}{(1+(n+1)|z|^2)^2}\\&&-\frac{3}{(1+|z|^2)^2} \nonumber
+2\frac{n+1}{n}        \left(\frac{n+1}{1+(n+1)|z|^2}-\frac{1}{1+|z|^2}
\right)dz\wedge                    d\zb\\                    \nonumber
&=&\frac{-4n}{2}+(n+1)-3+2\frac{n+1}{n}\log(n+1)
=-n-2+\left[2+\frac{2}{n}\right]\log(n+1).
\end{eqnarray*}
We are now able to compute
\begin{eqnarray*}
\widehat{F}_\star\left(        \widehat{c_1}(\overline{T\scn},\alpha_n)
\widehat{c_2}(\overline{T\scn},\alpha_n)\right)      &=&a\left(16+8\log
2\pi-4   \left(-1+\left[1+\frac{1}{n}\right]\log  (n+1)\right)  +4\log
2\pi\right.\\  &&\left.   +4\log 2\pi-\left(4n+4-\left[n+4+\frac{4}{n}
\right]\log(n+1)\right)   \right)\\  &=&a\left(n\log(n+1)-4n+16+16\log
2\pi\right).
\end{eqnarray*}

For analytic terms,
\begin{eqnarray*}
F_\star\left[Td(TS_n)   R(TS_n)ch(\mathcal{O})\right]  &=&\frac{2\zeta
'(-1)+\zeta   (-1)}{2}F_\star  c_1(TS_n)^2\\   &=&8\zeta  '(-1)+4\zeta
(-1))\\  F_\star\left[Td(TS_n)  R(TS_n)ch(\Omega_{S_n}^1)\right]&=&0\\
F_\star\left[Td(TS_n)    R(TS_n)ch(\Omega_{S_n}^2)\right]   &=&-8\zeta
'(-1)-4\zeta (-1)).
\end{eqnarray*}

\subsection{Conclusion}
We are now about to end our computations.
We just have to notice from application of the arithmetic
 Riemann-Roch theorem to the map $f : \p1\to \spec\zbb$ using
$\widehat{c_1}(\overline{T\p1},\omega_\pbu)=2\widehat{x}+a(\log
2\pi)$ that the analytic torsion of $\pbu$ is
$$\tau(\pbu,\omega_\pbu )=\frac{1+\log 2\pi}{3}-4\zeta '(-1)-2\zeta (-1).$$  
We find
\begin{theo}
\begin{eqnarray*}
\tau(\sn )
&=&\frac{n\log(n+1)}{24}-\frac{n}{6}
+\log\frac{n+2}{2}+2\tau(\pbu,\omega_\pbu).
\\
\tau(\sn,\Omega _{\sn}^1 )&=&0\\
\tau(\sn,\Omega _{\sn}^2 )&=&-\tau(\sn ).
\end{eqnarray*}
\end{theo}
Theses results are compatible with theorem 3.1 in ~\cite{RS}.

\section{Using Berthomieu-Bismut formula}
The second idea for the computation of the analytic torsion $\tau(S_n)$
 is to apply Berthomieu-Bismut formula~\cite{ber-bi} for the
composition of submersions~:
$$\sc_n=\pb (\mathcal{E}_n)\stackrel{\pi}\to\p1\stackrel{f}\to 
\spec\zbb.$$
We follow their construction in our particular setting.

\subsection{On Leray spectral sequence}
Leray spectral sequence of $\pi$ gives a canonical isomorphism $\sigma$
between the determinant line bundles. 
$$\lambda_f(R\pi_\star \osn)\to\lambda_F(\osn )$$
where we have set $F:=f\circ\pi$.
Note that we do not consider
their dual as in ~\cite{ber-bi}. 
Here, for $\pi$ is a locally trivial family of projective lines, 
$R\pi_\star \osn=\pi_\star \osn=\mathcal{O}_{\pb^1}$. The
spectral sequence hence degenerates in $E_2$ and the isomorphism
$\sigma$ is 
\begin{eqnarray*}
det H^0(\pb^1, \mathcal{O}_{\pb^1})
&\to & det H^0(S_n,\osn) \\
1'&\mapsto &1''=\pi^\star 1'
\end{eqnarray*}

\subsection{On Quillen metrics}
We now describe the Quillen metrics on determinant line bundles.
We choose the trivial metric on $\osn$.
The fibers of $\pi$ are endowed with the metric induced by 
$\alpha_n:=\Theta(\mathcal{O}_{E_n}(1),h_n)$.
 The bundle $R\pi_\star \osn=\mathcal{O}_{\pb^1}$ is endowed 
with its $L^2$ metric which is the trivial one for $\alpha_n$ is of
volume $1$ on every fiber of $\pi$.
Choose $\om_{\pb^1}:=x=\Theta(\mathcal{O}_\pbu (1),h)$
the Fubini-Study metric of volume $1$ as metric on $\pb^1$. 

Hence for $1'\in Rf_\star (R\pi_\star \osn)=H^0(\pb^1,\mathcal{O}_{\pb^1})$
\begin{eqnarray*}
\n 1'\n^2_{L^2}=\int_{\pb^1}1\om_{\pb^1}=1 &;&
\n 1'\n^2_Q=e^{-\tau(\pb^1)}.
\end{eqnarray*}
For $1''\in RF_\star \osn=H^0(\sn,\osn)$
\begin{eqnarray*}
\n 1''\n^2_{L^2}=\int_{S_n}1\frac{\alpha_n^2}{2}=\frac{n+2}{2}&;&
\n 1''\n^2_Q=\frac{n+2}{2}e^{-\tau(S_n)}.
\end{eqnarray*}
We can compute the Quillen norm of the isomorphism $\sigma$
\begin{eqnarray*}
\log\n\sigma\n^2_{\lambda^{-1 }_f(R\pi_\star \osn)\otimes\lambda_F(\osn )}
=\log\frac{\n 1''\n^2_Q}{\n 1'\n^2_Q}
=\tau(\pb^1)-\tau(S_n)+\log\frac{n+2}{2}.
\end{eqnarray*}

\subsection{Berthomieu-Bismut formula}
On the other hand this norm is computed by the Berthomieu-Bismut
formula
\begin{eqnarray*}
\log \n\sigma\n^2_{\lambda^{-1 }_f(R\pi_\star \osn)\otimes\lambda_F(\osn )} 
&=&-\int_{\pb^1} Td (T\pb^1,\om_{\pbu})Tors(\alpha_n,1)\\
&&+\int_{S_n}\widetilde{Td}(TS_n,T\pbu,\alpha_n,\om_\pbu)ch(\osn,1).
\end{eqnarray*}
Two kinds of secondary objects are used in this result. On one hand,
the form $Tors(\alpha_n,1)$ is the analytic torsion form of
Bismut-K{\"o}hler~\cite{bk} 
which fulfills the relation
\begin{eqnarray}
\label{torsion}&&
dd^c Tors(\alpha_n,1)
=\pi_\star(Td(T_{S_n/\pbu},\alpha_n)ch(\osn,1))-ch(R\pi_\star\osn,1).
\end{eqnarray}
On the other hand, 
the class $\widetilde{Td}(TS_n,T\pbu,\alpha_n,\om_\pbu)$ is the Bott-Chern
class for the Todd characteristic form ~\cite{gs-annals}
of the following metrized exact sequence

$$0\to (T_{S_n/\pbu},\alpha_n)\to (TS_n,\alpha_n)
\stackrel{d\pi}\to(\pi^\star T\pbu,\pi^\star\omega_\pbu)\to 0$$
which fulfills the relation
\begin{eqnarray*}
dd^c \widetilde{Td}(TS_n,T\pbu,\alpha_n,\om_\pbu)
=Td(T_{S_n/\pbu},\alpha_n)Td(\pi^\star T\pbu,\pi^\star\omega_\pbu)
-Td (TS_n,\alpha_n).
\end{eqnarray*}
Note that our sign conventions differ from the ones in~\cite{ber-bi} 
both in the definition of $\sigma$ and of the secondary objects.

\subsection{On the analytic torsion term}\label{todd}
The tool for the computation of the torsion form $Tors(\alpha_n,1)$ is
the arithmetic Riemann-Roch theorem gotten from equation
~(\ref{torsion}) by double transgression (see for example~\cite{B}
Theorem 4.4)
\begin{eqnarray}
\label{grr}\nonumber
a(Tors(\alpha_n,1))&=&
\widehat{\pi}_\star\left(\widehat{Td}^R(\overline{T_{\sc_n/\p1}})
\widehat{ch}(\overline{\oscn},1) \right)
-\widehat{ch}(\overline{R\pi_\star\oscn})\\
&=&\widehat{\pi}_\star\left(
  \widehat{Td}(\overline{T_{\scn/\p1}})\right)
-\widehat{\pi}_\star\left(\widehat{Td}(
\overline{T_{\scn/\p1}})a(R(T_{\sn/\pbu}))\right)-1.
\end{eqnarray}

We compute the arithmetic Todd class of $\overline{T_{\sc_n/\p1}}$
 using 
$\widehat{Td}=1+\frac{\widehat{c}_1}{2}+\frac{\widehat{c}_1^2}{12}$
for a line bundle on an arithmetic surface. From
 lemma~\ref{lem:classes}, we get
\begin{eqnarray*}
\widehat{Td}(\overline{T_{\scn/\p1}},\alpha_n )
&=&1+\widehat{\alpha_n}-\frac{n+2}{2}\widehat{x}+a(\frac{\log 2\pi}{2})
   \\
&&+\frac{1}{3}\widehat{\alpha_n}^2
-\frac{n+2}{3}\widehat{\alpha_n}\widehat{x}
+\frac{(n+2)^2}{12}\widehat{x}^2+\frac{\log 2\pi}{6}a(2\alpha_n-(n+2)x)\\
&=&1+\widehat{\alpha_n}-\frac{n+2}{2}\widehat{x}+a(\frac{\log 2\pi}{2})
   \\
&&+a\left(\frac{1+\log 2\pi}{3}\alpha_n
+\frac{n^2-2(n+2)\log 2\pi}{12}x
+\frac{1}{3}\frac{\langle \Theta (E^\star,h)a^\star,a^\star\rangle_h}
{\langle a^\star,a^\star\rangle_h}\right) .
\end{eqnarray*}
We compute the direct image,
$$\widehat{\pi}_\star(\widehat{Td}(\overline{T_{\scn/\p1}}))
=1+a(\frac{1+\log 2\pi}{3}).$$
The contribution of the $R$ class is purely analytic. 
First recall that $
c_1(T_{\scn/\p1},\alpha_n)=2\alpha_n-(n+2)\pi^\star x$ and
that $\pi_\star \alpha_n^2 =c_1(E_n,h_n)=(n+2)x$ so that 
$\pi_\star  c_1(T_{\scn/\p1},\alpha_n)^2=0$.
We find 
\begin{eqnarray*}
\lefteqn{\widehat{\pi}_\star(\widehat{Td}(\overline{T_{\scn/\p1}})
a(R(T_{\sn/\p1})))}\\ 
&=& a\left(
\pi_\star\left(Td(T_{S_n/\pbu},\alpha_n)R(T_{S_n/\pbu},\alpha_n)\right)
\right)\\ 
&=& a\left( \pi_\star\left( (1+\frac{1}{2}c_1( T_{S_n/\pbu}))
(2\zeta '(-1)+\zeta (-1))c_1( T_{S_n/\pbu})  \right)\right)\\
&=& a\left(4\zeta '(-1)+2\zeta (-1)\right).
\end{eqnarray*}
Back to formula ~(\ref{grr}), we find
$$Tors(\alpha_n,1)=\frac{1+\log 2\pi}{3}-4\zeta '(-1)-2\zeta (-1)
=\tau(\pbu,\omega_\pbu ).$$
We remark that the degree two part of the torsion vanishes.
We conclude this step
\begin{eqnarray*}
\int_{\pb^1} Td (T\pb^1,\om_{\pbu})Tors(\alpha_n,1)
&=& Tors(\alpha_n,1)\int_{\pb^1} Td (T\pb^1)\\
&=& Tors(\alpha_n,1)=\tau(\pbu ).
\end{eqnarray*}

\subsection{On the Bott-Chern term}
We now turn to the computation of 
$\widetilde{Td}(TS_n,T\pbu,\alpha_n,\om_\pbu)$.

For $ch(\osn,1)=1$, we only need to know 
$\widetilde{Td_3}(TS_n,T\pbu,\alpha_n,\om_\pbu)$. For
$Td_3=\frac{1}{24}c_1c_2$, we derive from ~\cite{gs-annals}
(prop.1.3.1.2) that
\begin{eqnarray}\label{eqn:todd}
\nonumber 24\widetilde{Td_3}(TS_n,T\pbu,\alpha_n,\om_\pbu)
&=&\widetilde{c_1}(TS_n,T\pbu,\alpha_n,\om_\pbu)
  c_2(T_{S_n/\pbu} \oplus \pi^\star T\pbu)\\
&&
+c_1(TS_n,\alpha_n ) \widetilde{c_2}(TS_n,T\pbu,\alpha_n,\om_\pbu).
\end{eqnarray}
The arithmetic relations between Chern classes in the following two
exact sequences 
$$
\begin{array}{ccccccccc}
 0 &\to&(T_{S_n/\pbu},\alpha_n)&\to&(TS_n,\alpha_n)&\to&
(\pi^\star T\pbu,\alpha_n) &\to &0 \\
0 &\to&(T_{S_n/\pbu},\alpha_n)&\to&(TS_n,\alpha_n)&\to&
(\pi^\star T\pbu,\pi^\star\omega_\pbu) &\to &0 \\
\end{array}$$ and the relation
\begin{eqnarray*}
\widehat{c_1}(\overline{\pi^\star T\p1},\alpha_n)
=\widehat{c_1}(\overline{\pi^\star T\p1},\pi^\star \omega_\pbu)
-a(\log\frac{{\alpha_n}_{|\pi^\star T\pbu}}{\pi^\star\omega_\pbu})
\end{eqnarray*}
enables us to infer
\begin{eqnarray*}
\widetilde{c}(TS_n,T\pbu,\alpha_n,\om_\pbu)
&=&\widetilde{c}(TS_n,T\pbu,\alpha_n,\alpha_n)
+c(T_{S_n/\pbu},\alpha_n)
\log\frac{{\alpha_n}_{|\pi^\star T\pbu}}{\pi^\star\omega_\pbu}.
\end{eqnarray*}

We can now evaluate the contribution of the Bott-Chern part. The first
term in formula~(\ref{eqn:todd}) is 
\begin{eqnarray*}
\lefteqn{\widetilde{c_1}(TS_n,T\pbu,\alpha_n,\om_\pbu)
  c_2(T_{S_n/\pbu} \oplus \pi^\star T\pbu)}&&\\
&=&\widetilde{c_1}(TS_n,T\pbu,\alpha_n,\om_\pbu) 
c_1(T_{S_n/\pbu},\alpha_n)
c_1(\pi^\star T\pbu,\pi^\star\om_\pbu )\\
&=&\log\frac{{\alpha_n}_{|\pi^\star T\pbu}}{\pi^\star\omega_\pbu}
\left(2\alpha_n-(n+2)x\right)2x
=4\log\frac{{\alpha_n}_{|\pi^\star T\pbu}}{\pi^\star\omega_\pbu}
\alpha_n x
\end{eqnarray*}
so that
\begin{eqnarray*}
\lefteqn{F_\star\left(\widetilde{c_1}(TS_n,T\pbu,\alpha_n,\om_\pbu)
  c_2(T_{S_n/\pbu} \oplus \pi^\star T\pbu)\right)}&&\\
&=&4\pi_\star\left(
\log\frac{1+(n+1)|z|^2\n e^\star\n^{2n}}{1+|z|^2\n e^\star\n^{2n}}
\alpha_n\right)f_\star\omega_{\pbu}\\
&=&4\int_\cb\log\frac{1+(n+1)|z|^2}{1+|z|^2}
\frac{i}{2\pi}\frac{dz\wedge d\zb}{(1+|z|^2)^2}
=\left[4+\frac{4}{n}\right]\log(n+1)-4.
\end{eqnarray*}
As for the second Bott-Chern class of $(TS_n,T\pbu,\alpha_n,\om_\pbu)$
\begin{eqnarray*}
 \widetilde{c_2}(TS_n,T\pbu,\alpha_n,\om_\pbu)
&=&\widetilde{c_2}(TS_n,T\pbu,\alpha_n,\alpha_n)
+ c_1(T_{S_n/\pbu},\alpha_n)
\log\frac{{\alpha_n}_{|\pi^\star T\pbu}}{\pi^\star\omega_\pbu}.
\end{eqnarray*}
This in turn cuts the computations into two pieces we already computed
in section ~\ref{RR}.
We find
\begin{eqnarray*}
F_\star\left(c_1(TS_n,\alpha_n )
  \widetilde{c_2}(TS_n,T\pbu,\alpha_n,\om_\pbu)\right) 
&=&4n+4-\left[n+4+\frac{4}{n} \right]\log(n+1)
\end{eqnarray*}
Summing up
\begin{eqnarray*}
\int_{S_n}\widetilde{Td}(TS_n,T\pbu,\alpha_n,\om_\pbu)ch(\osn,1)
&=&\frac{n}{6}-\frac{n\log(n+1)}{24}.
\end{eqnarray*}
This ends the second proof of our main formula.

\subsection{Remark on deformations}
The  only deformations  of the  surface $S_n$  are the  surfaces $S_m$
where  $m$ is  of the  same parity  as $n$.   Hirzebruch  surfaces are
pairwise diffeomorphic according to the parity of $n$, but the natural
diffeomorphisms as described in~\cite{MK} Part I, Theorem 4.2. 
are not isometric for  the metrics we choose which are
constructed  algebraically.  
More precisely, formulas for the deformation families of $S_{0}$ 
specializing to $S_{2p}$
$$
\begin{array}{ccc}
S_{0,t} &\stackrel{t\to 0}\longrightarrow &S_{2p}\\
\stackrel{\Psi_t}\simeq\downarrow & &\\
S_0&&
\end{array}
$$
lead, on appropriate open sets, to explicit formula for the metric on
$S_{0,t}$ transfered from the metric on $S_0$~:
$$\Psi_t^\star \alpha_0 =\pi^\star \omega_\pbu 
+ dd^c\log\left( 1+\left| \frac{x^p}{t}-\frac{1}{z}\right|^2\right)$$
and
$$\Psi_t^\star \alpha_0 =\pi^\star \omega_\pbu 
+ dd^c\log\left( 1+\left| \frac{z}{tx^pz+t^2}\right|^2\right).$$
This shows that the metric acquires singularities in the
specialization : it vanishes generically on generic fibers and tends
to Fubini-Study metric on the fiber over $0$.

Using  a genuine
diffeomorphism  between $S_{2p}$  and  $S_0$, the  computation of  the
analytic   torsion    of   $S_{2p}$   would    require   the   anomaly
formula~\cite{bgs} Theorem 1.23.

\section{Appendix}
We will use Hodge identity $\delta'' =[i\Lambda_{\alpha_n},d']$ to
check that the exhibited form is harmonic.

Starting from
\begin{eqnarray*}
\lefteqn{dd^c\log(1+|z|^2\n e^\star\n^{2n})}\\
&=&(n-\frac{n}{1+|z|^2\n e^\star\n^{2n}})\pi^\star \omega_\pbu
+\frac{\n e^\star\n^{2n}}{(1+|z|^2\n e^\star\n^{2n})^2}
\frac{i}{2\pi}\Big|dz+nzd'\log\n e^\star \n^2 \Big|^2
\end{eqnarray*}
and
\begin{eqnarray*}
\lefteqn{dd^c\log(1+(n+1)|z|^2\n e^\star\n^{2n})}\\
&=& (n-\frac{n}{1+(n+1)|z|^2\n e^\star\n^{2n}})\pi^\star \omega_\pbu
+\frac{(n+1)\n e^\star\n^{2n}}{(1+(n+1)|z|^2\n e^\star\n^{2n})^2}
\frac{i}{2\pi}\Big|dz+nzd'\log\n e^\star \n^2 \Big|^2
\end{eqnarray*}
we derive
\begin{eqnarray*}
\lefteqn{ dd^c\log\frac{{\alpha_n}_{|\pi^\star T\pbu}}{\pi^\star\omega_\pbu}
=dd^c\log\frac{1+(n+1)|z|^2\n e^\star\n^{2n}} {1+|z|^2\n e^\star\n^2}}\\
 &=&
\left(\frac{n}{1+|z|^2\n e^\star\n^{2n}}
-\frac{n}{1+(n+1)|z|^2\n e^\star\n^{2n}}\right)\pi^\star x\\
&&+ \left(\frac{(n+1)\n e^\star\n^{2n}}{(1+(n+1)|z|^2\n
    e^\star\n^{2n})^2}
-\frac{\n e^\star\n^{2n}}{(1+|z|^2\n e^\star\n^{2n})^2}\right)
\frac{i}{2\pi}\Big|dz+nzd'\log\n e^\star \n^2 \Big|^2\\
 &=&\frac{n^2|z|^2\n e^\star\n^{2n}}
{(1+|z|^2\n e^\star\n^{2n})(1+(n+1)|z|^2\n e^\star\n^{2n})}
\pi^\star \omega_\pbu\\
 &&+\frac{n\n e^\star\n^{2n}(1-(n+1)|z|^4\n e^\star\n^{4n})}
{(1+|z|^2\n e^\star\n^{2n})^2 (1+(n+1)|z|^2\n e^\star\n^{2n})^2}
\frac{i}{2\pi}\Big|dz+nzd'\log\n e^\star \n^2 \Big|^2.
\end{eqnarray*}
Recalling,
$$\alpha_n=\frac{1+(n+1)|z|^2\n e^\star\n^{2n}}
{1+|z|^2\n e^\star\n^{2n}}\pi^\star \omega_\pbu
+\frac{\n e^\star\n^{2n}}{(1+|z|^2\n e^\star\n^{2n})^2}
\frac{i}{2\pi}\Big|dz+nzd'\log\n e^\star \n^2 \Big|^2.$$
we are led to
\begin{eqnarray*}
\Lambda_{\alpha_n}\pi^\star \omega_\pbu &=& \frac{1+|z|^2\n e^\star\n^{2n}}
{1+(n+1)|z|^2\n e^\star\n^{2n}}
\end{eqnarray*}
and to
\begin{eqnarray*}
\Lambda_{\alpha_n} dd^c\log\frac{(1+(n+1)|z|^2\n e^\star\n^{2n})}
{1+|z|^2\n e^\star\n^{2n}} 
&=& n\frac{1-|z|^2\n e^\star\n^{2n}}{1+(n+1)|z|^2\n e^\star\n^{2n}}.
\end{eqnarray*}
This enables to compute
\begin{eqnarray*}
\Lambda_{\alpha_n}\left((n+2)\pi^\star\omega_\pbu
-dd^c\log\frac{(1+(n+1)|z|^2\n e^\star\n^{2n})}
{1+|z|^2\n e^\star\n^{2n}}\right) &=&2
\end{eqnarray*}
which proves that $\omega_H:=\pi^\star\omega_\pbu
-\frac{1}{n+2}dd^c\log\frac{\langle \Theta (E^\star,h)a^\star,a^\star\rangle}
{\pi^\star\omega_\pbu\langle a^\star,a^\star\rangle}$ is harmonic.
Remark that $\omega_H^2$ is $d$-exact and non-zero hence not harmonic
even if it is the product of two harmonic forms.

\end{document}